\documentclass[12pt]{amsart}
\usepackage{amssymb,latexsym,amsthm}
\usepackage{graphicx}

\newcommand{\inter}{\mathop{\mathrm{Int}}}

\newtheorem{theorem}{Theorem}[section]
\newtheorem{lemma}[theorem]{Lemma}
\newtheorem{cor}[theorem]{Corollary}

\theoremstyle{definition}

\newtheorem{example}[theorem]{Example}

\theoremstyle{remark}

\numberwithin{equation}{section}

\begin{document}

\baselineskip=17pt

\title[]
{A local Mazur-Ulam theorem}

\author{Osamu~Hatori}
\address{Department of Mathematics, Faculty of Science, 
Niigata University, Niigata 950-2181 Japan}
\curraddr{}
\email{hatori@math.sc.niigata-u.ac.jp}

\thanks{The author was partly 
supported by the Grants-in-Aid for Scientific 
Research, The 
Ministry of Education, Science, Sports and Culture, Japan.}

\keywords{normed spaces, isometries, the Mazu-Ulam theorem}


\maketitle

\begin{abstract}
We prove a local version of the Mazur-Ulam theorem. 
\end{abstract}
\section{Introduction}
In this paper we consider isometries between subsets of 
normed spaces. The Mazur-Ulam theorem asserts that 
an isometry from a normed space onto 
a normed space is real-linear up to translation (cf. 
\cite{mu,v}). An isometry from a normed space 
{\it into} a normed space need not be 
real-linear up to translation (cf. \cite{fi}). 
An isometry from an open set $U_1$ of a normed space 
onto an open set $U_2$ need not be extended to a 
real-linear map up to translation (see Example \ref{hoo}).
We show that if $U_1$ is star-shaped, then the isometry 
is extended to a real-linear map up to translation 
between the underlying normed spaces.
We also consider maps defined on a subset of a normed space 
which is not necessarily open.

Throughout the paper $B$ denotes a real normed space. For a 
subset $X$ of $B$, $\inter (X)$ is the interior of $X$.  
A star-shaped subset $K$ with a center $c$ of $B$ is
a set which satisfies that 
$tc+(1-t)x\in K$ for every $x\in K$ and $0\le t\le 1$.
Let $a\in B$ and $\varepsilon >0$. 
The open ball $\{x\in B:\|x-a\|<\varepsilon \}$ is 
denoted by $B_{\varepsilon}(a)$ and 
$\overline{B_{\varepsilon}(a)}$ its closure in $B$. 
For a pair $a$ and $b$ in $B$ the set 
$\{x\in B:ta+(1-t)b,\,0\le t\le 1\}$ is said to be a 
segment between $a$ and $b$ and is denoted by $[a,b]$.
\section{Isometries between open sets}
We begin by showing 
a preliminary lemma. We prove it by making use of an idea of 
V\"ais\"al\"a \cite{v}
\begin{lemma}\label{vvl}
Let $c\in B$ and a map $\psi:B\to B$ be defined as 
$\psi (z)=2c-z$ Suppose that $L$ is a non-empty 
bounded subset of $B$ such that $c\in L$ and $\psi (L)=L$. 
If ${\mathcal T}$is a 
surjective isometry from $L$ onto itself.
Then ${\mathcal T}(c)=c$.
\end{lemma}
\begin{proof}
Let $W$ be the set of all surjective isometries from $L$ onto 
itself. Note that $W$ is not empty since the identity function 
is in $W$. Let
\[
\lambda =\sup\{\|g(c)-c\|:g\in W\}.
\]
Since $L$ is bounded $\lambda<\infty$. 
We will show that $\lambda=0$. Suppose that 
$g\in W$. Let $g^*=g^{-1}\circ \psi \circ g$. Then 
$g^*\in W$. Hence 
\begin{multline*}
\lambda\ge \|g^*(c)-c\|=\|g^{-1}\circ \psi\circ g(c)
-c\|\\
=\|\psi \circ g(c)-g(c)\|=2\|g(c)-c\|.
\end{multline*}
It follows that $\lambda \ge 2\lambda$ since $g$ can be chosen 
arbitrary, hence $\lambda =0$.
\end{proof}
A real vector space with a metric 
$d(\cdot,\cdot)$ satisfying $d(a+u,b+u)=d(a,b)$ 
for all $a,b,u$, and for which addition and 
scalar multiplication are jointly continuous is 
called a metric real vector space.
\begin{lemma}\label{lmu}
Let ${\mathcal B}_1$ be a real normed space and 
${\mathcal B}_2$ a metric real vector space with a 
metric $d(\cdot,\cdot)$. Suppose that 
$U_1$ and $U_2$ are non-empty open subsets of 
${\mathcal B}_1$ and 
${\mathcal B}_2$ respectively. Suppose that ${\mathcal T}$ is 
a surjective isometry 
\rm{(}$d({\mathcal T}(a),{\mathcal T}(b))=\|a-b\|$ for 
every $a,b\in U_1$\rm{)} from $U_1$ onto $U_2$ 
and $f,g\in U_1$. If $f$ and $g$  satisfy 
the equation $(1-r)f+rg\in U_1$ for every 
$r$ with $0\le r \le 1$, then the equality
\[
{\mathcal T}(\frac{f+g}{2})=\frac{{\mathcal T}(f)+{\mathcal T}(g)}{2}
\]
holds.
\end{lemma}
\begin{proof}
Let $h,h'\in U_1$. Suppose that there exists $\varepsilon >0$ 
which satisfies that 
$\frac{\|h-h'\|}{2}<\varepsilon$, and 
\[
\{a\in {\mathcal B}_1:\|a-h\|<\varepsilon,\,\,\|a-h'\|<\varepsilon \}\subset U_1,
\]
\[
\{u\in {\mathcal B}_2:d(u,{\mathcal T}(h))<\varepsilon,\,\,
d(u,{\mathcal T}(h'))\|<\varepsilon \}\subset U_2.
\]
We will show that 
${\mathcal T}(\frac{h+h'}{2})=\frac{{\mathcal T}(h)+{\mathcal T}(h')}{2}$.
Set 
$r=\frac{\|h-h'\|}{2}$ and let
\[
L_1=
\{a\in {\mathcal B}_1:\|a-h\|=r=\|a-h'\|\},
\]
\[
L_2=
\{u\in {\mathcal B}_2:d(u,{\mathcal T}(h))=r=
d(u,{\mathcal T}(h'))\}.
\]
Set also $c_1=\frac{h+h'}{2}$ and 
$c_2=\frac{{\mathcal T}(h)+
{\mathcal T}(h')}{2}$. Then we have ${\mathcal T}(L_1)=L_2$, 
$c_1\in L_1 \subset U_1$, and $L_2 \subset U_2$. 
Let 
\[
\psi_2(y)=2c_2-y \quad (y\in {\mathcal B}_2).
\]
Then $\psi_2$ is an isometry and $\psi_2(L_2)=L_2$.

Let $Q={\mathcal T}^{-1}\circ\psi_2\circ{\mathcal T}$. 
Then $Q$ is well-defined and is a surjective isometry 
from  $L_1$ onto itself. 
Then by Lemma \ref{vvl} $Q(c_1)=c_1$. Henceforth 
${\mathcal T}(c_1)=c_2$.

We assume that $f$ and $g$ are as described. Let
\[
K=\{(1-r)f+rg:0\le r\le 1\}.
\]
Since $K$ and ${\mathcal T}(K)$ are compact, there is $\varepsilon >0$ with
\[
\inf\{\|a-b\|:a\in K,\,b\in {\mathcal B}_1\setminus U_1\}
>\varepsilon, 
\]
\[
\inf\{d(u,v):u\in {\mathcal T}(K),\,
v\in {\mathcal B}_2\setminus U_2\}>
\varepsilon.
\]
Then for every $h\in K$ we have 
\[
\{a\in {\mathcal B}_1:\|a-h\|<\varepsilon\}\subset U_1
\]
and
\[
\{u\in {\mathcal B}_2:\|u-{\mathcal T}(h)\|
<\varepsilon \}\subset U_2.
\]
Choose a natural number $n$ with $\frac{\|f-g\|}{2^n}<\varepsilon$. 
Let 
\[
h_k=\frac{k}{2^n}(g-f)+f
\]
for each $0\le k\le 2^n$. By the first part of the proof we have
\[
{\mathcal T}(h_k)+{\mathcal T}(h_{k+2})-2{\mathcal T}(h_{k+1})=0
\qquad \text{($k$)}
\]
holds for $0\le k\le 2^n-2$. For $0\le k\le 2^n-4$,
adding the equations ($k$), 2 times of ($k+1$), and ($k+2$) we have
\[
{\mathcal T}(h_k)+{\mathcal T}(h_{k+4})-2{\mathcal T}(h_{k+2})=0,
\]
whence the equality 
\[
{\mathcal T}(\frac{f+g}{2})=\frac{{\mathcal T}(f)+{\mathcal T}(g)}{2}
\]
holds by induction on $n$.
\end{proof}
An isometry between open sets of normed spaces need not 
be extended to a linear map up to translation.
\begin{example}\label{hoo}
Let ${\mathcal X}
=\{x,y\}$ be a compact Hausdorff space consisting of two points 
and $C({\mathcal X})$ denote the Banach algebra of all 
complex-valued continuous functions on ${\mathcal X}$.
Let 
\[
{\mathcal U}=\{f\in C({\mathcal X}):\|f\|<1\}\cup
\{f\in C(X):\|f-f_0\|<1\},
\]
where $f_0\in C({\mathcal X})$ is defined as $f_0(x)=0,\,f_0(y)=10$.
Suppose that 
\[
{\mathcal T}:{\mathcal U}\to {\mathcal U}
\]
is defined as ${\mathcal T}(f)= f$ if $\|f\|<1$ and 
${\mathcal T}(f)=\tilde f$ if $\|f-f_0\|<1$, where 
\begin{equation*}
\tilde f(t)=
\begin{cases}
-f(t),& t=x \\
f(t),& t=y.
\end{cases}
\end{equation*}
Then ${\mathcal T}$ is an isometry from ${\mathcal U}$ onto itself, 
while it cannot be extended to a real linear isometry up to translation.
\end{example}
\section{Sufficient conditions for real-linearity}
\begin{theorem}\label{ss}
Let $B_1$ and $B_2$ be real normed 
spaces. Let $U_1$ be a non-empty star-shaped open subset of 
$B_1$. Suppose that $T:U_1\to B_2$ is an isometry such that 
$T(U_1)$ is open in $B_2$. 
Then there exists a surjective real-linear isometry 
$\tilde T_0$ from $B_1$ onto $B_2$ and 
$u\in B_2$ such that 
\[
T(a)=\tilde T_0(a)+u\quad (a\in U)
\]
holds.
\end{theorem}
\begin{proof}
Let $a_0$ be a center of $U_1$; i.e., $[a_0,x]\subset U_1$ for 
every $x\in U_1$. 
Let $\varphi:B_1\to B_1$ be defined as $\varphi_1(x) =x+a_0$ 
($x\in B_1$) and 
$\varphi_2:B_2\to B_2$ defined as 
$\varphi_2(x)=x-T(a_0)$ ($x\in B_2$). Then 
$\varphi_2\circ T\circ \varphi_1:U_1-a_0\to B_2$ is an 
isometry such that 
$\varphi_2\circ T\circ \varphi_1(U_1-a_0)$ is open in $B_2$. 
We will show that 
$\varphi_2\circ T\circ \varphi_1$ is extended to a surjective 
real-linear isometry $\tilde T_0$ from $B_1$ onto $B_2$. 
It will follow that the conclusion holds. 
Let $V_1=U_1-a$ and $T_0=
\varphi_2\circ T\circ \varphi_1$. 
There is $r>0$ with 
$B_{3r}(0)\subset V_1$. 
Let $a\in V_1$. Since $0$ is a center of 
$V_1$, $ta=ta+(1-t)0\in V_1$. Then by 
Lemma \ref{lmu}
\[
T_0\left(\frac{a}{2}\right)=
T_0\left(\frac{a+0}{2}\right)=
\frac{T_0(a)+T_0(0)}{2}=\frac{T_0(a)}{2}
\]
holds. For every 
$0\le s_1\le 1$ and $0\le s_2\le 1$, 
$t(s_1a)+(1-t)(s_2a)\in V_1$, hence
\[
T_0\left(
\frac{s_1a+s_2a}{2}\right)=
\frac{T_0(s_1a)+T_0(s_2a)}{2}
\]
holds. Applying the above equation and by induction on $n$ 
\[
T_0\left(\frac {m}{2^n}a\right)=\frac{m}{2^n}T_0(a)
\]
holds for every positive integer $n$ and 
$m=0,1,\dots , 2^n$. It follows that 
\begin{equation}\label{eq1}
T_0(ta)=tT_0(a)
\end{equation}
holds for every $a\in V_1$ and $0\le t\le 1$. 
Let $a,b\in \overline{B_r(0)}$. Then 
$a+b\in B_{3r}(0)\subset V_1$ holds, and 
$ta+(1-t)b\in \overline{B_r(0)}\subset V_1$ ($0\le t\le 1$) 
and (\ref{eq1}) imply the 
equations 
\begin{equation}\label{eq2}
T_0(a+b)=2T_0\left(\frac{a+b}{2}\right)=
2\times \frac{T_0(a)+T_0(b)}{2}=
T_0(a)+T_0(b)
\end{equation}
by Lemma \ref{lmu}. Define a map $\tilde T_0:B_1\to B_2$ as 
\begin{equation*}
\tilde T_0(x)= 
\begin{cases}
0,& x=0, \\
\frac{\|x\|}rT_0\left(\frac{r}{\|x\|}x\right), & x\ne 0.
\end{cases}
\end{equation*}
We will show that $T_0(x)=\tilde T_0(x)$ for 
every $x\in V_1$. Let $x\in V_1$. If $x=0$, then 
$T_0(x)=0=\tilde T_0(x)$. 
Suppose that $r\le \|x\|$. Then by (\ref{eq1})
\[
T_0\left(\frac{r}{\|x\|}x\right)=
\frac{r}{\|x\|}T_0(x),
\]
so $\tilde T_0(x)=T(x)$ holds. Suppose that 
$x\ne 0$ and $\|x\|<r$. Then 
\begin{equation}\label{3.25}
T_0(x)=T_0\left(\frac{\|x\|}{r}\frac{r}{\|x\|}x\right)
=\frac{\|x\|}{r}T_0\left(\frac{r}{\|x\|}x\right)=
\tilde T_0(x)
\end{equation}
hold by (\ref{eq1}).

We will show that 
\begin{equation}\label{eq3.5}
\tilde T_0(sx)=s\tilde T_0(x)
\end{equation}
for every $x\in B_1$ and $s\in {\mathbb R}$. 
Let $a\in \overline{B_r(0)}$. 
Then $-a\in \overline{B_r(0)}$, and by (\ref{eq2}) 
$T_0(-a)=-T_0(a)$ holds. Let 
$x\in B_1$ and $s\in {\mathbb R}$. 
If $x=0$ or $s=0$, then 
$\tilde T_0(sx)=s\tilde T_0(x)$ holds. 
Suppose that $x\ne 0$ and $s>0$. 
Then 
\[
\tilde T_0(sx)=
\frac{\|sx\|}{r}T_0\left(
\frac{rs}{\|sx\|}x\right)=
\frac{s\|x\|}{r}T_0\left(\frac{r}{\|x\|}x\right)
=s\tilde T_0(x).
\]
Suppose that $x\ne 0$ and $s<0$. Then 
\[
\tilde T_0(sx)=\frac{-s\|x\|}{r}
T_0\left(
\frac{-r}{\|x\|}x
\right)
=\frac{s\|x\|}{r}T_0\left(
\frac{r}{\|x\|}x\right)=
s\tilde T_0(x)
\]
since $\frac{r}{\|x\|}x\in \overline
{B_r(0)}$.

We will show that 
\begin{equation}\label{eq4}
\tilde T_0(x+y)=
\tilde T_0(x)+\tilde T_0(y)
\end{equation}
holds for every $x,y\in B_1$. If 
$x+y=0$, then $y=-x$ and hence (\ref{eq4}) holds by 
(\ref{eq3.5}). 
Suppose that $x+y\ne 0$. If $x=0$ or $y=0$, then 
(\ref{eq4}) holds since $\tilde T_0 (0)=0$. 
Suppose that $x\ne 0$ and $y\ne 0$. Then by (\ref{eq2}) 
and (\ref{3.25}) 
\begin{multline*}
\tilde T_0\left(
\frac{r}{\|x\|+\|y\|}x+\frac{r}{\|x\|+\|y\|}y\right)
=
T_0\left(
\frac{r}{\|x\|+\|y\|}x+\frac{r}{\|x\|+\|y\|}y\right) \\
=
T_0\left(\frac{r}{\|x\|+\|y\|}x\right)+
T_0\left(\frac{r}{\|x\|+\|y\|}y\right)\\
=
\tilde T_0\left(\frac{r}{\|x\|+\|y\|}x\right)+
\tilde T_0\left(\frac{r}{\|x\|+\|y\|}y\right)
\end{multline*}
follows and (\ref{eq4}) holds by (\ref{eq3.5}).

We will show that $\tilde T_0$ is a surjection. 
Let $y\in B_2$. Then there is 
$r_0>0$ such that 
$\frac{y}{r_0}\in T_0(V_1)$ since $T_0(V_1)$ is open 
and $0=T_0(0)\in T_0(V_1)$. 
Then there is $x_0\in V_1$ with 
$T_0(x_0)=\frac{y}{r_0}$. It follows that 
$\tilde T_0(r_0x_0)=y$. 

We will show that $\tilde T_0$ is an isometry. Let 
$x\in B_1$. If $x=0$, then 
$\|\tilde T_0(x)\|=\|x\|$. Suppose that $x\ne 0$. Then 
\[
\tilde T_0(x)=\frac{\|x\|}{r}\tilde T_0
\left(\frac{r}{\|x\|}x\right)=
\frac{\|x\|}{r}T_0\left(
\frac{r}{\|x\|}x\right)
\]
hence
\begin{equation*}
\|\tilde T_0(x)\|=
\frac{\|x\|}{r}\|T_0\left(\frac{r}{\|x\|}x\right)\|\\
=\frac{\|x\|}{r}\|\frac{r}{\|x\|}x-0\|=\|x-0\|=\|x\|
\end{equation*}
since $T_0$ is an isometry and $T_0=\tilde T_0$ on 
$V_1$. Thus $\tilde T_0$ is 
an isometry since $\tilde T_0$ is linear.
\end{proof}
There are two preliminary lemmata for the following corollary. 
They may be standard, 
but proofs are included for the sake of 
completeness.
\begin{lemma}\label{lem1}
Suppose that $U$ is a non-empty open subset of $B$ 
and $p\in B\setminus U$. Then 
\[
V_{p,U}=\cup_{x\in U}[p,x] \setminus \{p\}
\]
is an open subset of $B$.
\end{lemma}
\begin{proof}
Suppose that $x_0\in V_{p,U}$. Then there exist 
$y_0\in U$ and $0\le t\le 1$ with 
$x_0=tp+(1-t)y_0$. Note that $t<1$ holds, in fact, for 
$p \not\in V_{p,U}$. Since $U$ is open there exists 
$\varepsilon >0$ with $B_{\varepsilon}(y_0)\subset U$. 
Then by a simple calculation 
\[
B_{(1-t)\varepsilon}(x_0)\subset V_{p,U}
\]
holds. Hence $V_{p,U}$ is open for $x_0$ is arbitrary.
\end{proof}
\begin{lemma}\label{lem2}
Suppose that $X$ is a convex subset of $B$ and $\inter (X)
\ne \emptyset$. Then $\inter (X)$ is also convex and the 
closure $\overline{\inter (X)}$ of $\inter (X)$ contains 
$X$.
\end{lemma}
\begin{proof}
We will show that $[a,b]\subset \inter (X)$ for every pair 
$a$ and $b$ in $\inter (X)$. 
Suppose that $a,b \in \inter (X)$. Then there exists 
an open neighbourhood $U\subset \inter(X)$ of $b$ with 
$a\not\in U$. 
Then Lemma \ref{lem1} insures that 
$V_{a,U}$ is open. 
Since $X$ is convex, $[a,x]\subset X$ for every $x\in U$, 
hence $V_{a,U}\subset X$, hence 
$V_{a,U}\subset \inter(X)$ since $V_{a,U}$ is open. 
It follows that $[a,b]\subset \inter(X)$ since $a\in 
\inter (X)$; $\inter(X)$ is convex. 

Let $p\in X$. 
Let $a$ be any element in $\inter(X)$. Then there 
exists an open neighbourhood $U$ of $a$ such that 
$p\not\in U\subset \inter(X)$. 
Hence $V_{p,U}\subset \inter(X)$ since 
$V_{p,U}\subset X$ for $X$ is convex and 
$V_{p,U}$ is open
by Lemma \ref{lem1}. 
Let $x_n=(1-\frac{1}{n})p+\frac{1}{n}a$ for each positive 
integer $n$. Then 
$x_n\in V_{p,U}$ and $x_n\to p$ as $n\to \infty$. Then 
$p\in \overline{\inter(X)}$.
\end{proof}
\begin{cor}
Let $B_1$ and $B_2$ be real normed spaces. Let $X$ be a 
convex subset of $B_1$ and $\inter(X)\ne \emptyset$. 
Suppose that $T:X\to B_2$ is isometric and 
$T(\inter(X))$ is open in $B_2$. 
Then $T$ is extended to a real-linear isometry up to 
translation.
\end{cor}
\begin{proof}
Since the restriction $T|_{\inter(X)}:\inter(X)\to 
B_2$ of $T$ to $\inter(X)$ satisfies the hypotheses of 
Theorem \ref{ss}, $T|_{\inter(X)}$ is extended to a surjective 
real-linear isometry up to translation 
$\widetilde{T|_{\inter(X)}}$ from $B_1$ onto $B_2$. 
Since $T$ and $\widetilde{T|_{\inter(X)}}$ are isometric and 
$X\subset \overline{\inter(X)}$ holds by Lemma \ref{lem2}, 
$\widetilde{T|_{\inter(X)}}=T$ on $X$.
\end{proof}
An isometry from a star-shaped closed subset need not be 
extended to a real-linear map up to translation.
\begin{example}
Let ${\mathbb R}^2_{\max}={\mathbb R}^2$ 
as real linear spaces and the norm is defined as 
$\|(x,y)\|=\max \{|x|,|y|\}$ for $(x,y)\in 
{\mathbb R}^2_{\max}$.
Let 
\[
X_1=\{(x,y)\in 
{\mathbb R}^2_{\max}:-1\le x\le 0,\,|y|\le -x\},
\]
\[
X_2=\{(x,0)\in {\mathbb R}^2_{\max}:0\le x\le 1\}
\]
and
\[
X=X_1\cup X_2.
\]
Then $X$ is a star-shaped closed subset of 
${\mathbb R}^2_{\max}$. Let 
\[
T:X\to {\mathbb R}^2_{\max}
\]
be defined as 
\begin{equation*}
T((x,y))=
\begin{cases}
(x,y), & (x,y)\in X_1, \\
(x,\sin x ),& (x,y)\in X_2.
\end{cases}
\end{equation*}
Then $T$ is isometry and $T(\inter(X))$ is open in 
${\mathbb R}^2_{\max}$. 
On the other hand $T$ is not extended to a linear map since
\[
T(-(1,0))=(-1,0)\ne (-1,-\sin 1)=-T((1,0))
\]
and $T((0,0))=(0,0)$.
\end{example}



\end{document}